\documentclass{IEEEtran4PSCC}

\usepackage{amsmath,amsthm}
\usepackage{amsfonts}
\usepackage{amssymb}

\usepackage{tabularx}
\usepackage{graphicx}
\usepackage{booktabs}
\usepackage{longtable}
\usepackage{lscape}
\usepackage{tabu}
\usepackage{color}
\usepackage{pst-all, pst-node}

\usepackage[ruled,vlined,linesnumbered]{algorithm2e}

\SetAlgorithmName{Algorithm Schema}{\autoref}{List of Algorithms}

\newcommand{\stt}{\text{ s.t. }}

\newcommand\Sos{{\operatorname{\text{\small $\Sigma$}}}}

\usepackage{ulem}
\let\underbar\uline

\newcommand{\A}{\mathcal{A}} 
\newcommand{\N}{\mathcal{N}} 
\newcommand{\K}{\mathcal{K}} 
\newcommand{\J}{\mathcal{J}} %
\renewcommand{\L}{\mathcal{L}} %
\newcommand{\M}{\mathcal{M}} %

\newcommand{\Y}{\mathcal{Y}} %

\newcommand{\iunit}{\mathbf{j}} 

\begin{document}
\title{MINLP in Transmission Expansion Planning }

\author{
\IEEEauthorblockN{Jakub Mare{\v c}ek \\ Martin Mevissen \\ Jonas Christoffer Villumsen}
\IEEEauthorblockA{IBM Research -- Ireland \\ Dublin, Ireland}
}	

\maketitle

\begin{abstract}
Transmission expansion planning requires forecasts of 
demand for electric power and a model of the underlying physics, i.e., power flows. 
We present three approaches to deriving exact solutions to the transmission expansion planning problem in the alternating-current model,
for a given load.
\end{abstract}

\begin{IEEEkeywords}
Power system analysis computing, Optimization, Numerical Analysis (Mathematical programming)
\end{IEEEkeywords}

\section{Introduction}

Consider the problem of optimal investment in line capacity, such that 
a sum of annualised investment costs and an estimate of operational costs is minimised. 
The estimate of operational costs and its computational complexity
depend on the model of power flows used.
Considering the recent progress in the development of convergent solvers for 
polynomial optimisation \cite{Ghaddar15}, 
we explore the options for solving the transmission expansion problem exactly
in the alternating-current model.

This is motivated by the observation \cite{Zhang2013} that the quality of the approximation of the alternating-current model has a major impact on the investment decisions. 
Specifically, the use of the simplistic direct-current approximation (DCOPF)
may result in no lines being built. 
Various piece-wise linearisations may result in various lines being built, 
 other than those built considering the alternating-current model (ACOPF) proper. 
This is the case even when loads are known exactly, 
 i.e., independently of the uncertainty in the load.

We compare three convergent approaches to solving the transmission expansion planning 
problem in the alternating-current model.
First, we study both the current-voltage (IV) and power-voltage (PQV) formulations of the problem 
  as polynomial optimisation problems and derive semidefinite-programming relaxations (SDP) thereof
  using the techniques of \cite{Ghaddar15}. 
Second, we introduce a novel lift-and-branch-and-bound procedure using SDP relaxations we introduce,
  which makes it possible to obtain global optima for small instances of the transmission expansion problem.
Finally, we compare these approaches with state-of-the-art piece-wise linearisations based on the current-voltage formulation and a rudimentary DCOPF approximation. 
Although we do not consider multiple scenarios for the demand, it would be easy to extend the
work in the direction of two-stage or multi-stage stochastic programming.


\section{The Problem}

Formally, let introduce the problem using:
\begin{itemize}
\item $(\N,\A)$ be the graph representing an electrical network with buses $\N$ and lines $\A$ 
\item $p_i + \iunit q_i$ be the complex net power injection at bus $i \in \N $,
\item $z_i + \iunit w_i$ be the complex net current injected at bus $i \in \N $,
\item $z_{ij} + \iunit w_{ij}$ be the complex current flow on line $(i,j) \in \A $ (with slight abuse of notation), and let
\item $v_i + \iunit u_i$ be the complex voltage at bus $i \in \N$ and 
\item $o_{ij}=1$ if circuit $(i,j) \in \A $ is open and 0 otherwise.
\end{itemize}

The transmission expansion problem may be stated as follows:
\begin{align}
        \min c^{\top}p & & \tag*{[IV]} \label{eq:ACOPFIV-0}\\
\stt    & z_i + \sum_{j\in \N} (z_{ji} - z_{ij}) = 0 & \forall i \in \N \label{eq:ACOPFIV-1}\\
        & w_i + \sum_{j\in \N} (w_{ji} - w_{ij}) = 0 & \forall i \in \N \label{eq:ACOPFIV-2} \\
        & \underbar{p}_i  \leq p_i  = v_iz_i + u_iw_i \leq \bar{p}_i & \forall i \in \N \label{eq:ACOPFIV-5}\\
        & \underbar{q}_i  \leq q_i  = u_iz_i - v_iw_i \leq \bar{q}_i & \forall i \in \N \label{eq:ACOPFIV-6} \\
        & \underbar{v}_i^2  \leq v_i^2 + u_i^2 \leq \bar{v}_i^2 & \forall i \in \N \label{eq:ACOPFIV-7}\\
        & o_{ji} \in \{0,1\} & \forall (i,j) \in A \label{ACOPFIV-10}
\end{align}
\begin{align}
        o_{ij}=0 \Rightarrow & v_i - v_j = R_{ij}z_{ij} - X_{ij}w_{ij} & \forall (i,j) \in \A \label{eq:ACOPFIV-3}\\
        o_{ij}=0 \Rightarrow & u_i - u_j = X_{ij}z_{ij} + R_{ij}w_{ij} & \forall (i,j) \in \A \label{eq:ACOPFIV-4} \\
        o_{ij} = 0 \Rightarrow & z_{ij}^2 + w_{ij}^2 \leq \bar{z}_{ij}^2 & \forall (i,j) \in \A \label{eq:ACOPFIV-8i} \\
        o_{ij} = 1 \Rightarrow & z_{ij}^2 + w_{ij}^2 \leq 0 & \forall (i,j) \in \A \label{eq:ACOPFIV-8j}
\end{align}

where:
\begin{itemize}
\item $c$ denotes the vector of generator marginal costs, which are assumed to be known ans constant, for simplicity,  
\item $R_{ij}$ and $X_{ij}$ denote the resistance and reactance respectively of line $(i,j)$.
\item $\underbar{p}$ and $\bar{p}$ denote minimum active power generation and capacity vectors, 
\item $\underbar{q}$ and $\bar{q}$ denote minimum reactive power generation and capacity vectors, 
\item $\underbar{v}$ and $\bar{v}$ denote vectors of minimum and maximum voltage magnitudes, and 
\item $\bar{z}$ is the thermal capacity limit of lines.
\end{itemize}


\section{A Piece-Wise Linearisation}

Several authors have proposed piece-wise linearisations of the current-voltage formulation
of optimal power flows (ACOPF-IV) and related problems,
which yield mixed-integer linear programming formulations \cite{Alguacil2003,Fisher2008,Oneill2012, Ferreira2013,Trodden2014}. 
The idea of piece-wise linearisation is well-known 
and often very efficient.
Most of the proposed piece-wise linearisations, e.g. \cite{Fisher2008,Zhang2013,Trodden2014},
however, cannot reach the global optimum in the limit of the number of 
segments, due to the additional assumptions surveyed in Table~1 of \cite{Trodden2014}. 
In the following, we extend the formulation proposed in \cite{Ferreira2013} to 
accommodate dispatch of real and reactive power, and
hence obtain a ``principled'' piece-wise linearisation.

The upper bound on voltage magnitude in \eqref{eq:ACOPFIV-7} are convex quadratic
constraints, that may be replaced by an outer linear approximation as in \cite{Oneill2012}.
The remaining non-linear constraints \eqref{eq:ACOPFIV-5}, \eqref{eq:ACOPFIV-6}, and the
first inequality of \eqref{eq:ACOPFIV-7} is discussed below.

In general, the complex current injection may be written in terms of voltage and power as

$$
I = z + \iunit w = \frac{S}{V} = \frac{p + \iunit q}{v + \iunit u} = \frac{vp+uq}{v^2+u^2} + \iunit \frac{vq-up}{v^2+u^2}
$$

By employing a discretisation of the complex two-dimensional $(V, S)$-space, we can 
linearise the complex current injection. 
We propose to discretise the complex voltage space $V$ along its polar
coordinates, while discretising the power along rectangular coordinates.
This ensures a tight approximation to the feasible area of the voltage
space.
Let $\hat{|V|}$, $\hat{\theta}$, $\hat{p}$, and $\hat{q}$ be
the vectors of values at the discretisation points in the four real dimensions, 
voltage magnitude, voltage angle, real, and reactive power, respectively.
We then evaluate the value of the current injection in each of the 
$|\hat{|V|}|\times|\hat{\theta}|\times|\hat{p}|\times |\hat{q}|$ discretisation points as 
$\hat{z} + \iunit \hat{w}$.

Now, at any point in the $(|V|,\theta,p,q)$-space we can approximate all relevant quantities
as a convex combination $\lambda$ of the values in the immediately surrounding discretisation points.
That is, for all $i \in \N$

\begin{align*}
        v_i &= \sum_{j,k,l,m}\hat{|V|}_i^{j} \cos \hat{\theta}_i^k \lambda_i^{jklm} &
        u_i &= \sum_{j,k,l,m}\hat{|V|}_i^{j} \sin \hat{\theta}_i^k \lambda_i^{jklm} \\ 
        p_i &= \sum_{j,k,l,m}\hat{p}_i^{l} \lambda_i^{jklm} &
        q_i &= \sum_{j,k,l,m}\hat{q}_i^{m} \lambda_i^{jklm} \\
        z_i &= \sum_{j,k,l,m} \hat{z}_i^{jklm} \lambda_i^{jklm} &      
        w_i &= \sum_{j,k,l,m} \hat{w}_i^{jklm} \lambda_i^{jklm} \\
        |V_i| &= \sum_{j,k,l,m} \hat{|V|}_i^{j} \lambda_i^{jklm} & &
\end{align*}

In each dimension we choose a convex combination of discretisation points,

\begin{align*}
        \sum_{j\in\J} \sum_{k\in\K} \sum_{l\in\L} \sum_{m\in\M} \lambda_i^{jklm} = 1  \quad \forall i \in \N \\
\end{align*}

with $\lambda_i^{jklm}  \geq 0$ for all $i \in \N, j\in\J, k\in\K, l\in\L, m\in\M$.
The so called ``special ordered set of type 2'' (SOS2) constraints ensures that we choose a 
combination of the two closest discretisation points in each dimension:

\begin{align*}
        \lambda^{jklm} & \leq \psi^{j-1} + \psi^j, \forall j\in \J \setminus \{0\} \\
        \lambda^{jklm} & \leq \chi^{k-1} + \chi^k, \forall k\in \K \setminus \{0\} \\
        \lambda^{jklm} & \leq \varphi^{l-1} + \varphi^l, \forall l\in \L \setminus \{0\} \\
        \lambda^{jklm} & \leq \upsilon^{m-1} + \upsilon^m, \forall m \in \M \setminus \{0\} \\
        \lambda^{j,0,0,0} &\leq \psi^0 & 
        \lambda^{0,k,0,0} &\leq \chi^0 \\
        \lambda^{0,0,l,0} &\leq \varphi^0 &
        \lambda^{0,0,0,m} &\leq \upsilon^0 \\
        e^{\top}\psi &= e^{\top}\chi = e^{\top}\varphi= e^{\top}\upsilon = 1 & 
        \psi, \chi, \varphi, \upsilon &\in \{0,1\} \\
\end{align*}

where $\J, \K, \L$, and $\M$ are the index sets of the discretisation points $\hat{|V|}, \hat{\theta}$, and $\hat{p}$
$\hat{q}$, respectively. That is,
\begin{itemize}
\item $\psi^j=1$ if and only if $V$ is in the interval $\left[ \hat{|V|}^{j},\hat{|V|}^{j+1}\right]$,
\item $\chi^k=1$ if and only if $\theta$ is in the interval $\left[ \hat{\theta}^{k},\hat{\theta}^{k+1}\right]$, 
\item $\varphi^l=1$ if and only if $p$ is in the interval $\left[ \hat{p}^{l},\hat{p}^{l+1}\right]$, while
\item $\upsilon^m=1$ if and only if $q$ is in the interval $\left[ \hat{q}^{m},\hat{q}^{m+1}\right]$.
\end{itemize}

For the reference bus, the voltage is fixed and only the two dimensional real
$(p,q)$-space is discretised, while for demand buses without generation,
 power is fixed and only the voltage space is discretised.

We can model line-use by introducing the binary variable $o_{ij}=1$ if and only if the line $(i,j)$ is not available, i.e.~not installed or the switch is open. That is, we replace  the implications (\ref{eq:ACOPFIV-3}--\ref{eq:ACOPFIV-8j}) by the disjunctive constraints:

\begin{align}
        -M o_{ij} \leq R_{ij}z_{ij} - X_{ij}w_{ij} + v_j - v_i  \leq M o_{ij} \quad  \forall (i,j) \in \A \\
        -M o_{ij} \leq X_{ij}z_{ij} + R_{ij}w_{ij} + u_j - u_i  \leq M o_{ij} \quad \forall (i,j) \in \A \\
        z_{ij}^2 + w_{ij}^2 \leq \bar{z}_{ij}^2 (1 - o_{ij}) \quad  \forall (i,j) \; \in \A \label{eq:ACOPFIV}
\end{align}

In summary, a ``principled'' piece-wise linearisation
 converges to the true optimum in the large limit of the number of 
 evaluation points, in theory.
This comes at the price of discretisation of the two-dimensional $(V, S)$-space, where we discretise the complex voltage space $V$ along its polar
coordinates, while discretising the power along rectangular coordinates,
which seems superior to the alternative choices of discretisation.

\section{Sum-of-Squares Approaches}

Alternatively, one can exploit a rich history of research into polynomial optimisation: 
\begin{align}
\min \quad& f(x) \notag \\
\mbox{s.t.  }& g_i(x) \geq 0 \qquad i=\{1,\dots,m\} \tag*{[PP]}
\end{align}
as surveyed in \cite{Handbook} and elsewhere. 
Let us use $\mathcal{P}_d(S)$ to denote the cone of polynomials of degree at most $d$ that are non-negative over some
$S\subseteq \mathbb{R}^n.$
A homogeneous polynomial $h(x)$ of degree $2d$ in $n$-dimensional vector $x$ is sum-of-squares (SOS, \cite{Choi1995}) if and only if there exist homogeneous polynomials of degree $d$, denoted $g_1(x),\ldots,g_k(x)$ such that 
$h(x)=\sum_{i=1}^k g_i(x)^2$. 
for conductance, as usual in the power systems community.)
We use $\Sigma_d$ to denote the cone of polynomials of degree at most $d$ that are sum-of-squares of polynomials.
It has been shown \cite{Choi1995} that each $\mathcal{P}_{2d}(\mathbb{R}^n)$ can be approximated as closely as desired by a sum-of-squares of polynomials,
in the $l_1$-norm of its coefficient vector, albeit with a possibly large $k$.
Using $\mathcal{G}=\{g_i(x): i=1,\dots,m \}$ and denoting $S_{\mathcal{G}}=\{x \in \mathbb{R}^n :  g(x) \geq 0, \; \forall g \in \mathcal{G}\}$ the basic closed semi-algebraic set defined by $\mathcal{G}$, Lasserre reformulates [PP] as
\begin{align}
\max \quad & \varphi & \mbox{s.t.  }& f(x)-\varphi  \geq 0 \quad \forall \: x \in S_{\mathcal{G}}, \notag \\
= \max \quad & \varphi & \mbox{s.t.  } &f(x)-\varphi \in \mathcal{P}_d(S_{\mathcal{G}}). \tag*{[PP-D]}
\end{align}
which allows for approximations up to arbitrary accuracy within a number of hierarchies \cite{Lasserre1,Lasserre2006}.

The so called ``dense'' hierarchy of Lasserre \cite{Lasserre1} approximates $\mathcal{P}_d(S_{\mathcal{G}})$
by the cone $\mathcal{{K}}^r_{\mathcal{G}}$, where
\begin{align}
\mathcal{{K}}^{r}_{\mathcal{G}} =  \Sos_{r}+\sum_{i=1}^{m}g_i(x) \Sos_{r-\deg(g_i)}, \label{eqn:K} 
\end{align}
and $r \geq d$. 
The corresponding optimisation problem over $S$ can be written as:
{\begin{align}
\max_{\varphi, \sigma_i(x)} \:  & \varphi \tag*{[PP-H$_r$]$^*$} \label{eq-Lass} \\
\mbox{s.t. } & f(x)-\varphi= \sigma_0(x)+\sum_{i=1}^{m} \sigma_i(x)g_i(x) \notag\\
&\sigma_0(x) \in \Sos_{r}, \: \sigma_i(x) \in \Sos_{r-\deg(g_i)}. \notag  
\end{align}}
and \ref{eq-Lass} can be reformulated as a semidefinite optimisation problem. We denote the dual of [PP-H$_r$]$^*$ by [PP-H$_r$], in keeping with 
previous work \cite{Ghaddar15}.

The so called ``sparse'' hierarchy of Waki et al. \cite{WKKM,Lasserre2006,KojimaMuramatsu2009} is based on the correlative sparsity of a polynomial optimisation problem [PP] of dimension $n$, which can be represented by the $n\times n$ correlative sparsity pattern matrix:
$$\mathcal{R}_{ij} = \begin{cases} 
\star& \text{for } i=j\\ 
\star& \text{for } x_i, x_j \text{ in the same monomial of } f \\ 
\star& \text{for } x_i, x_j \text { in the same constraint } g_k \\ 
0 & \text {otherwise},
\label{eq:CSP}
\end{cases} $$
and its associated adjacency graph $G$, the {\it correlative sparsity pattern graph}. Let $\{I_k \}_{k=1}^p$ be the set of maximal cliques of a chordal extension of $G$ following the construction in \cite{WKKM}, i.e. $I_k \subset \{ 1,\ldots,n \}$. 
The sparse approximation of $\mathcal{P}_d(S)$ is $\mathcal{{K}}^{r}_{\mathcal{G}}(I)$, given by
\begin{align*}
  \mathcal{{K}}^{r}_{\mathcal{G}}(I) = \sum_{k=1}^p \left(\Sigma_r(I_k) +  \sum_{j\in J_k} g_j \Sigma_{r-\deg(g_j)}(I_k)\right),
\end{align*}
where $\Sigma_d(I_k)$ is the set of all sum-of-squares polynomials of degree up to $d$ supported on $I_k$ and $(J_1,\ldots, J_p)$ is a partitioning of the set of polynomials $\{g_j\}_j$ defining $S$ such that for every $j$ in $J_k$, the corresponding $g_j$ is supported on $I_k$. The support $ I \subset \{1,\ldots,n\} $ of a polynomial contains the indices $i$ of terms $x_i$ which occur in one of the monomials of the polynomial.  The sparse hierarchy of SDP relaxations is then given by 
\begin{align}
&\max_{\varphi,\sigma_k(x), \sigma_{r,k}(x)} \:   \varphi \tag*{[PP-SH$_r$]$^*$} \label{eq-Lass-sparse}  \\ 
&\mbox{s.t. } f(x)-\varphi=  \sum_{k=1}^p \left(\sigma_{k}(x) +  \sum_{j\in J_k} g_j(x) \sigma_{j,k}(x)\right) \notag \\
&\quad \sigma_{k}\in\Sos_r((I_k)), \sigma_{j, k} \in \Sos_{r-\deg(g_j)}(I_k). \notag 
\end{align}
We denote the dual of [PP-SH$_r$]$^*$ by [PP-SH$_r$], again in keeping with 
previous work \cite{Ghaddar15}.


One can easily see the current-voltage (IV) formulation (\ref{eq:ACOPFIV-0}--\ref{eq:ACOPFIV-8j}) as a polynomial optimisation problem: 
one only needs to reformulate the constraints on line-use variable $o_{ij}$ to $o_{ij}^2 = o_{ij}$ and the implication with 
antecedent $o_{ij} = 0$ using either ``Big M'' constraints or perspective reformulation \cite{gunluk2010perspective}.
Subsequently, one can derive two hierarchies of semidefinite programming relaxations, as described above.
While these relaxations of the degree-2 polynomial optimisation problem are tractable, 
they also turn out to be rather weak.


One may also consider the power-voltage (PQV) formulation, where variables are complex voltages $\mathcal{V}$ at each bus and power at each generator. 
With matrices $Y_{k}, Y_{ij}$, $\bar{Y}_{k}, \bar{Y}_{ij}$ derived from the admittance matrix in the usual 
fashion \cite{Ghaddar15}, and without the line-use decision, 
the PQV formulation of optimal power flows is a polynomial optimisation problem of degree 2 or 4.
The line-use decision $o_{ij}=0$ in the antecedent of the implications of the investment problem, 
however, requires replacing constant matrices with variables such as, 
\begin{align}
\Y_{ij}(o) = \begin{cases}  y_{ii} + \sum_{i \neq j} {o_{ij} y_{ij}}, & \mbox{if } i = j \\  
-y_{ij} o_{ij},  & \mbox{if } i \neq j, (i, j) \in \A \\
0,  & \mbox{if } i \neq j, (i, j) \not\in \A \\
\end{cases},
\end{align}
For the thermal limits, one would hence have to introduce:
\begin{align}
(\text{tr}(\Y(o)Y_{ij}\mathcal{V}\mathcal{V}^T))^2+(\text{tr}(\Y(o)\bar{Y}_{ij}\mathcal{V}\mathcal{V}^T))^2 \notag \\
\leq \bar{z}_{ij}^2 (1 - o_{ij}) &\quad  \forall (i,j) \in A
\end{align}
which raises the degree of the polynomial optimisation problem to 9. 
This produces a strong relaxation,
albeit hard to use with general-purpose polynomial optimisation techniques.

\section{Lift-and-Branch-and-Bound}


Finally, we present a lift-and-branch-and-bound scheme.
Similarly to branch-and-bound approaches \cite{land1960automatic} in mixed-integer linear programming (MILP), 
we consider repeatedly a subproblem, where values of certain variables $o_{ij}$ 
are fixed to either 0 or 1. 
We denote the set of constraints fixing these variables to the prescribed values by $E$.
Outside of variables fixed in $E$, constraints $o_{ij}^2 = o_{ij}$ are relaxed to $0 \le o_{ij} \le 1$ in either the IV or PQV formulation above,
which are made progressively tighter in a newly introduced outer loop.
We denote such a sub-problem by [Relaxation-SH$_r(E)$], where Relaxation is either PQV or IV and $r$ is the counter of the outer loop,
i.e., the order of the relaxation.

\begin{algorithm}[tb!]
\caption{Lift-and-Branch-and-Bound} 
\label{alg:BranchLift}
\DontPrintSemicolon
\SetKwInOut{Input}{input}
\SetKwInOut{Output}{output}
initialise the best known upper bound $z \longleftarrow \infty $ \label{step:initz} \;
initialise a queue $Q \longleftarrow Q_{\text{init}}$, e.g. $Q_{\text{init}} = \emptyset$ \label{step:initQ} \;
initialise a queue of $Q' \longleftarrow \emptyset$ \;
initialise Relaxation to either PQV or IV \;
initialise $r$ to the minimum order required by Relaxation \label{step:initr} \;
\While{ $Q \cup Q' \not = \emptyset $  }{ \label{step:outerwhile}
 \While{ $Q \not = \emptyset $  }{ \label{step:innerwhilewhile}
    $E \longleftarrow  Q$.pop()\;
    \eIf{ $|E| = |\A|$, i.e. we have a complete solution in terms of $o_{ij}$ }{
      \eIf{ [PP$_2$-SH$_r$] of \cite{Ghaddar15} is feasible with cost $\underbar z$ } {
          \eIf{ $\underbar z \ge z$ } {
            Drop $E$
           } {
             \eIf{ flat-extension conditions \cite{curto1998flat} are satisfied } {
              Record the solution, update $z$ if needed, and drop $E$
             } {
              $Q'$.push($E$)
             }
           }
      } {
          Drop $E$
      }
    } {
      \eIf{ [Relaxation-SH$_r(E)$] is feasible with cost $\underbar z$ } { \label{step:secondif}
        \eIf{ solution to [Relaxation-SH$_r(E)$] has $o_{ij} \in \{0, 1\} \; \forall (i,j) \in A$  } {
             \eIf{ flat-extension conditions \cite{curto1998flat} are satisfied } {
              Record the solution, update $z$ if needed, and drop $E$
             } {
              $Q'$.push($E$)
             }
        } {
          \eIf{ $\underbar z \ge z$ } {
            Drop $E$
           } {
            Pick $(i,j) \in A$ not in $E$ such that the value of $o_{ij}$ in the solution is fractional \;
            $Q$.push($E \cup \{ o_{ij} = 0 \}$) \;
            $Q$.push($E \cup \{ o_{ij} = 1 \}$)
           }
        }
      } {
          Drop $E$
      }
    }
  }
  $Q \longleftarrow Q'$ \;
  $Q' \longleftarrow \emptyset $ \;  
  $r \longleftarrow r + 1$ \;
}
Return $z$ and the corresponding solution \;
\end{algorithm}

The pseudo-code of the algorithm is displayed in Algorithm Schema~\ref{alg:BranchLift}.
Initially, one considers the so called ``root relaxation'', where $E = \emptyset$.
Counter $r$ is initialised to the minimum required by either [PQV-SH$_r(E)$] or [IV-SH$_r(E)$].
A queue $Q$ stores subproblems, which define a partial or complete solution in terms of the investment decision.
While processing $E \in Q$, we may arrive at one of the four outcomes:
\begin{itemize}
\item a feasible solution is found 
\item infeasibility or a bound sufficiently strong to prune $E$ is found,    
\item if there is a variable $o_{ij}$, whose value is not fixed, $Q$ is extended with branches $o_{ij} = 0/1$ 
\item processing is deferred to $Q'$.
\end{itemize}

Specifically, if the set of constraints define a complete solution
in terms of the investment decision $o_{ij}$, the processing considers [PP$_2$-SH$_r$] of Ghaddar et al.\ \cite{Ghaddar15}.
Alternatively, starting on Line \ref{step:secondif}, we consider [Relaxation-SH$_r(E)$].
Once there are no elements left in $Q$, we move the contents of $Q'$ into $Q$, increment $r$, and repeat.
The hope is that we prune as many subproblems as possible using low-rank relaxations, so as to process fewer nodes for higher values of counter $r$.



\section{Computational Illustrations}
\label{sec:results}

We have evaluated the approaches on variants of a simple two-bus instance \cite{Bukhsh2013} and variants of Garver's \cite{Garver1970}
six-bus network. 
Throughout, e.g. in Table~\ref{tab:1res} 
below, 
we detail:
\begin{itemize} 
\item Ops: the per-hour costs of operations in US dollars
\item Obj: the sum of the per-hour costs of operations with the per-hour amortisation of the investment in US dollars
\item T: run-time of the solver in seconds, as measured on 
a machine equipped with 80 cores of Intel Xeon CPU E7 8850 and circa 700 GB of RAM
\item RMSE: the root mean squared error for the voltage magnitude VM$_i'$ and voltage angle VA$_i' \forall i \in \N$ of the solution in questions,
      compared to the voltage magnitude VM$_i$ a voltage angle VA$_i \forall i \in \N$ at a global optimum with the same fixed phase, 
     \begin{align}
     \sqrt{\frac{\sum_{i=1}^{|\N|} (\textrm{VM}_i - \textrm{VM}_i')^2 + (\textrm{VA}_i - \textrm{VA}_i')^2 }{2 |\N|}}
     \end{align}
\item MAPE: mean absolute percentage error of the solution in terms of voltage magnitude VM$_i'$ and voltage angle VA$_i' \forall i \in \N$ of the solution in questions,
      compared to the voltage magnitude VM$_i$ a voltage angle VA$_i \forall i \in \N $at a global optimum with the same fixed phase, 
      \begin{align}\frac{1}{2|\N|}\sum_{i=1}^{|\N|}  \left|\frac{\textrm{VM}_i' - \textrm{VM}_i}{\textrm{VM}_i}\right| + \left|\frac{\textrm{VA}_i' - \textrm{VA}_i}{\textrm{VA}_i}\right|
      \end{align}
\item B$i$: voltage magnitude (VM) and angle (VA) at the respective bus.
\end{itemize}
of the following solvers:
\begin{itemize} 
\item DC: the direct-current optimal power flow, as implemented in Matpower \cite{matpower}
\item PWL1: a piece-wise linearised model \ref{eq:ACOPFIV-0} with 
a coarse discretisation  of Table~\ref{tab:discretisations}
\item PWL2: a piece-wise linearised model \ref{eq:ACOPFIV-0} 
with a fine discretisation of Table~\ref{tab:discretisations}
\item IPM: an efficient interior-point method, as implemented in Matpower \cite{matpower}
\item SDP: the PQV SDP relaxation \cite{molzahn2011,lavaei2012zero}.
\end{itemize}

Wherever applicable, the infeasibility has been checked by both Matpower, the interior-point method and the {\tt insolvablepfsos\_limitQ} routine of Mohlzahn in Matpower 5.0. 
The interior-point method provides a heuristic, but widely-used indication thereof, 
 whereas the latter certifies the same.
Wherever applicable, global optimality has been certified by the rank of the solution of the SDP relaxation, as per \cite{lavaei2012zero,molzahn2011}.


\subsection{Case2}

Our computational experiments start with a 2-bus instance by Bukhsh et al. \cite{Bukhsh2013}, 
where one can invest into one line, with $r = 0.04$ and $x = 0.2$.
Trivially, there is 1 feasible configuration. 
First, we observe that both the Lavaei-Low relaxation (for ACOPF) and the lowest orders of both IV and PQV relaxations (for TEP) fail to find the global optimum, although PQV at $r = 5$ gets close.
Second, considering a degree-9 polynomial is involved in the PQV formulation, 
solving the relaxation with $r = 5$ requires circa 20 GB of RAM.
Further, we observe that the use of $o_{ij}^2 = o_{ij}$ and similar 
does not improve the relaxation considerably.
Even this trivial instance hence suggests that 
there are substantial limitations to the performance of the relaxations.

\begin{table}[t!]
\centering
\begin{tabular}{|l|r|r|r|r|}\hline
Piece-wise linearisation & Va (per 90 deg) &	Vm	&	P	&	Q	\\
\hline
crude	&	3	&	1	&	1	&	1	\\
fine	&	6	&	2	&	2	&	2	\\
very fine&	1000	&	2	&	2	&	2	
\\ \hline 
\end{tabular}
\caption{The discretisations used in the three piece-wise linearisations for case2 and case2mod.}
\label{tab:discretisations}
\end{table}


We have also introduced a 2-bus test instance (case2mod), 
where there are no existing lines, 2 parallel lines one can invest into, 
and hence 4 possible configurations. 
The two parallel lines between buses 1 and 2 differ in their admittance.
The first one has $r = 0.04$ and $x = 0.2$, while the second one has 
$r = 0.02$ and $x = 0.1$.
Both crude and fine piece-wise linearisations,
the interior point method implemented in Matpower \cite{matpower}, and
the SDP relaxation of Lavaei and Low \cite{molzahn2011,lavaei2012zero}
fail to find the global optimum.
When one adds the investment decision, the 
relaxation grow very quickly. 

Specifically, the size of the constraint matrix for case2mod as an investment problem grows from 
$513 \times 3311$ with 3449 non-zeros for the first order, which can be solved in 7.5 seconds,
to $16905 \times 163506$ (183646 non-zeros) for the second order, which requires 4655.8 seconds,
and beyond. It is not known how high in the hierarchy one would need to go to obtain the exact solution
to the investment problem.
Compare this to 
$138 \times 331$ matrix with 624 for the usual PQV SDP relaxation of Lavaei-Low \cite{lavaei2012zero,molzahn2011},
which can be solved in 1.65 seconds, but remains a challenge to extend to the investment decision.
We will provide the details in an extended version of the paper.


\begin{table*}
\centering
Buses:\\[2mm]
\begin{tabular}{|l|l|r|r|r|r|r|r|r|r|r|r|r|}\hline
bus & type &   $P_d$   &   $Q_d$  &  $G_s$  &  $B_s$  & area &  $V_m$  &  $V_a$  & baseKV & zone & $V_{\max}$ & $V_{\min}$ \\ \hline 
1 & 3 &  80.00 & 16.00 & 0.00 & 0.00 & 1.00 & 1.00 & 0.00 & 230.00 & 1.00 & 1.05 & 0.95 \\ \hline 
2 & 1 & 240.00 & 48.00 & 0.00 & 0.00 & 1.00 & 0.00 & 0.00 & 230.00 & 1.00 & 1.05 & 0.95 \\ \hline 
3 & 2 &  40.00 &  8.00 & 0.00 & 0.00 & 1.00 & 1.00 & 0.00 & 230.00 & 1.00 & 1.05 & 0.95 \\ \hline 
4 & 1 & 160.00 & 32.00 & 0.00 & 0.00 & 1.00 & 0.00 & 0.00 & 230.00 & 1.00 & 1.05 & 0.95 \\ \hline 
5 & 1 & 240.00 & 48.00 & 0.00 & 0.00 & 1.00 & 1.00 & 0.00 & 230.00 & 1.00 & 1.05 & 0.95 \\ \hline 
6 & 2 &  0.00  &  0.00 & 0.00 & 0.00 & 1.00 & 0.00 & 0.00 & 230.00 & 1.00 & 1.05 & 0.95 \\ \hline 
\end{tabular}\\[6mm]
Generators:\\[2mm]
\begin{tabular}{|c|c|c|c|c|c|c|c|c|c|c|c|c|c|}\hline
bus &  Qmax  &  Qmin  &  Vg  &  Pmax  & Pmin &  Pc1 &  Pc2 & Qc1min & Qc1max & Qc2min & Qc2max \\ \hline 
1 & 48.25 & -10.00 & 1.00 & 160.00 & 0.00 & 0.00 & 0.00 &  0.00  &  0.00  &  0.00  &  0.00  \\ \hline 
3 & 101.25 & -10.00 & 1.00 & 370.00 & 0.00 & 0.00 & 0.00 &  0.00  &  0.00  &  0.00  &  0.00  \\ \hline 
6 & 183.00 & -10.00 & 1.00 & 610.00 & 0.00 & 0.00 & 0.00 &  0.00  &  0.00  &  0.00  &  0.00  \\ \hline 
\end{tabular}\\[6mm]
Branches:\\[2mm]
\begin{tabular}{|l|l|r|r|r|r|r|r|r|}\hline
fbus & tbus & $R$ & $X$ & rateA  & rateB & rateC &	ratio & angle \\  \hline 
1 & 2 &   0.040 &  0.40 &   180 &  250 & 250 & 0 & 0 \\ \hline 
1 & 4 &   0.060 &  0.60 &   150 &  250 & 250 & 0 & 0 \\ \hline 
1 & 5 &   0.010 &  0.10 &   360 &  250 & 250 & 0 & 0 \\ \hline 
2 & 3 &   0.020 &  0.20 &   180 &  250 & 250 & 0 & 0 \\ \hline 
2 & 4 &   0.040 &  0.40 &   180 &  250 & 250 & 0 & 0 \\ \hline 
2 & 6 &   0.015 &  0.15 &   360 &  250 & 250 & 0 & 0 \\ \hline 
3 & 5 &   0.010 &  0.10 &   360 &  250 & 250 & 0 & 0 \\ \hline 
3 & 6 &   0.024 &  0.24 &   360 &  250 & 250 & 0 & 0 \\ \hline 
4 & 6 &   0.008 &  0.08 &   360 &  250 & 250 & 0 & 0 \\ \hline
\end{tabular}
\caption{The details of the instance depicted in Figure~\ref{fig:garver6instance} in the Matpower format. Columns not listed are uniformly at the default values.}
\label{tab:1}
\end{table*}

\subsection{Garver6y}

\begin{figure}[tb!]
\includegraphics[width=0.49\textwidth,clip=true,trim=4.5cm 20cm 4cm 4.5cm]{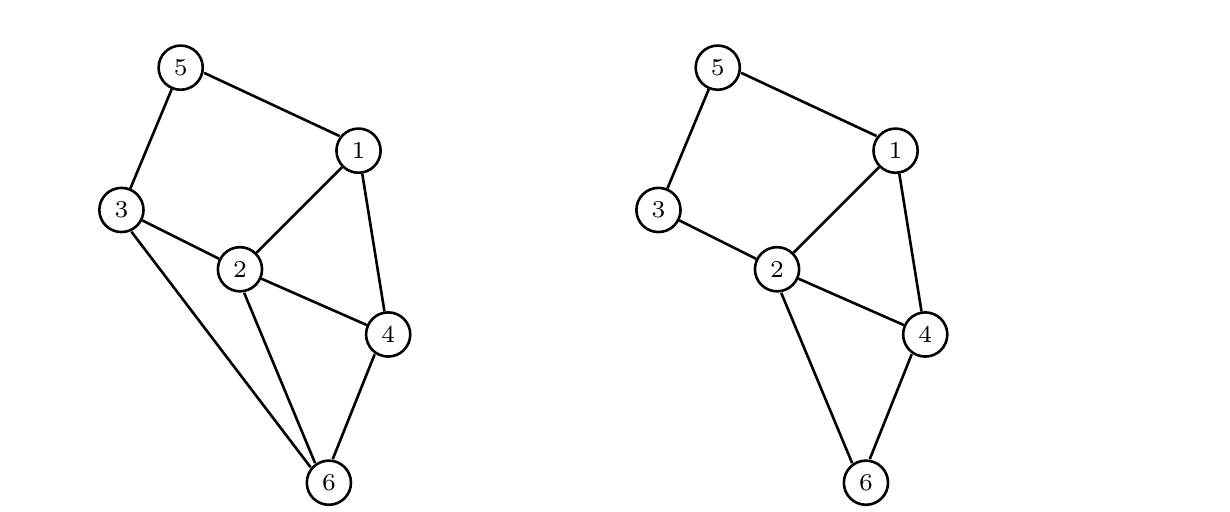}
\caption{The instance Garver6y with the AC-optimal solution (left) and 
the solution to PWL1 (right).
}
\label{fig:garver6instance}
\end{figure}

Next, we have 
introduced a small test instance based on the one by Garver \cite{Garver1970}.
As in Garver's original example we 
consider the connection of bus 6 to the existing system. 
For the sake of clarity, we consider three double circuit 
lines 2-6, 3-6, and 4-6, and hence 8 
possible configurations. 
The existing lines are
complemented by an extra circuit along line 1-5 and 3-5 and 
the corresponding lines are replaced by their equivalent 
single circuit line. 
The 
network configuration is shown in Figure \ref{fig:garver6instance}. 
Line investment costs, loosely based on 
amortisation to a per-hour cost assuming a 40-year planning
period and an interest rate of 3~\%,
are 100, 80, and 50 for the lines 2--6, 3--6, and 4--6, respectively.
See Table~\ref{tab:1} for details.

Table~\ref{tab:1conf} provides an overview of the 8 configurations.
We label the configurations 000 to 111, where the first-listed binary digit
indicates whether line $2 \to 6$ is built, the second-listed binary digit indicates
whether the line $3 \to 6$ is built, and the third-listed binary digit indicates 
whether $4 \to 6$ is built.
Configurations 000 and 010 are both AC and DC infeasible.
A DC model considers configuration 001 feasible with cost 1840,
100 feasible with cost 1840, 
110 feasible with cost 1412,
and 011 feasible with cost 1360.
Neither of those is AC feasible, though.
The remaining two configurations are AC feasible,
as detailed in Table~\ref{tab:1res},
but the piece-wise linearisation gives rather different solution from the
exact optimum recovered by the SDP relaxation.
Notably, the DC model would result in a different investment decision
from the crude piece-wise linearisation, 
which would be different still from the decision made using either 
the interior-point method or the semidefinite programming,
which both produce the globally optimal solution, in this case.

\begin{table}[t!]
\begin{tabular}{|c|c|c|c|r|r|r|r|r|}\hline
   Id & 2-6 & 3-6 & 4-6 & DC            & PWL1          & PWL2        & IPM        & SDP     \\ \hline 
   000           & 0         & 0         & 0         &               &               &             &            &         \\ \hline 
   001           & 0         & 0         & 1         & 1890       &               &             &            &         \\ \hline 
   010           & 0         & 1         & 0         &               &               &             &            &         \\ \hline 
   011           & 0         & 1         & 1         & {\bf 1490}       &               &             &            &         \\ \hline 
   100           & 1         & 0         & 0         & 1940       &               &             &            &         \\ \hline 
   101           & 1         & 0         & 1         & 1510       & {\bf 1936}    & 1868     & 1887  & 1887 \\ \hline 
   110           & 1         & 1         & 0         & 1592       &               &             &            &         \\ \hline 
   111           & 1         & 1         & 1         & 1590       & 1949       & {\bf 1770}     & {\bf 1818 }  & {\bf 1818} \\ \hline 
\end{tabular}
\caption{The 8 possible configurations of the system, depending on the investment into lines, and their respective operating costs
in dollars per hour,
whenever feasible within the model.}
\label{tab:1conf}
\end{table}

\begin{table*}[tb]
\centering
Configuration 101: \\
\begin{tabular}{|c|c|c|c|c|c|c|c|c|c|c|}\hline
& \multicolumn{2}{c|}{DC} & \multicolumn{2}{c|}{PWL1} & \multicolumn{2}{c|}{PWL2} & \multicolumn{2}{c|}{IPM} & \multicolumn{2}{c|}{SDP} \\
             &  VM  &   VA  & VM &   VA  & VM &   VA  &  VM &   VA  &  VM &   VA  \\ \hline 
   Obj / T   & 1360.00 &  0.04 & 1786.11 &  0.18 & 1717.66 &  1.39 & 1736.84 &  0.16 & 1736.84 &  1.04 \\ \hline 
 RMSE / MAPE &   5.91  &  1.21 &   2.63  &  0.67 &   0.08  &  0.02 &   0.00  &  0.00 &   0.00  &  0.00 \\ \hline 
      B1     &   1.00  &  0.00 &   0.96  &  0.00 &   0.98  &  0.00 &   0.98  &  0.00 &   0.98  &  0.00 \\ \hline 
      B2     &   1.00  &  6.75 &   0.93  &  5.12 &   0.95  &  1.06 &   0.95  &  1.21 &   0.95  &  1.21 \\ \hline 
      B3     &   1.00  & -5.79 &   1.01  &  4.19 &   1.02  &  0.69 &   1.01  &  0.80 &   1.01  &  0.80 \\ \hline 
      B4     &   1.00  & 21.00 &   0.96  & 14.57 &   0.98  &  9.64 &   0.98  &  9.76 &   0.98  &  9.76 \\ \hline 
      B5     &   1.00  & -9.77 &   0.94  & -4.92 &   0.95  & -6.63 &   0.95  & -6.74 &   0.95  & -6.74 \\ \hline 
      B6     &   1.00  & 33.98 &   1.03  & 24.94 &   1.04  & 19.38 &   1.05  & 19.51 &   1.05  & 19.51 \\ \hline 
\end{tabular}
 \\[6mm]
Configuration 111: \\
\begin{tabular}{|c|c|c|c|c|c|c|c|c|c|c|}\hline
& \multicolumn{2}{c|}{DC} & \multicolumn{2}{c|}{PWL1} & \multicolumn{2}{c|}{PWL2} & \multicolumn{2}{c|}{IPM} & \multicolumn{2}{c|}{SDP} \\
             &  VM  &   VA  & VM &   VA  & VM &   VA  &  VM &   VA  &  VM &   VA  \\ \hline 
   Obj / T   & 1360.00 &  0.03 & 1719.35 &  0.23 & 1540.48 &  2.84 & 1587.86 &  0.09 & 1587.86 &  0.88 \\ \hline 
 RMSE / MAPE &   1.75  &  2.01 &   1.74  &  2.21 &   0.29  &  0.17 &   0.00  &  0.00 &   0.00  &  0.00 \\ \hline 
      B1     &   1.00  &  0.00 &   0.98  &  0.00 &   0.98  &  0.00 &   0.98  &  0.00 &   0.98  & -0.00 \\ \hline 
      B2     &   1.00  &  1.72 &   0.94  & -1.77 &   0.95  &  0.20 &   0.96  &  0.07 &   0.96  &  0.07 \\ \hline 
      B3     &   1.00  &  0.26 &   1.02  &  3.90 &   1.01  &  1.79 &   1.01  &  2.24 &   1.01  &  2.24 \\ \hline 
      B4     &   1.00  & 10.38 &   0.96  &  3.68 &   0.98  &  7.79 &   0.99  &  7.23 &   0.99  &  7.23 \\ \hline 
      B5     &   1.00  & -6.74 &   0.96  & -4.89 &   0.95  & -6.17 &   0.95  & -6.00 &   0.95  & -6.00 \\ \hline 
      B6     &   1.00  & 20.83 &   1.02  & 12.35 &   1.05  & 17.04 &   1.05  & 16.39 &   1.05  & 16.39 \\ \hline 
\end{tabular}

\caption{The two AC-feasible configurations of Figure~\ref{fig:garver6instance}.
}
\label{tab:1res}
\end{table*}

\begin{figure}[tb!]
\centering
\includegraphics[width=0.45\textwidth]{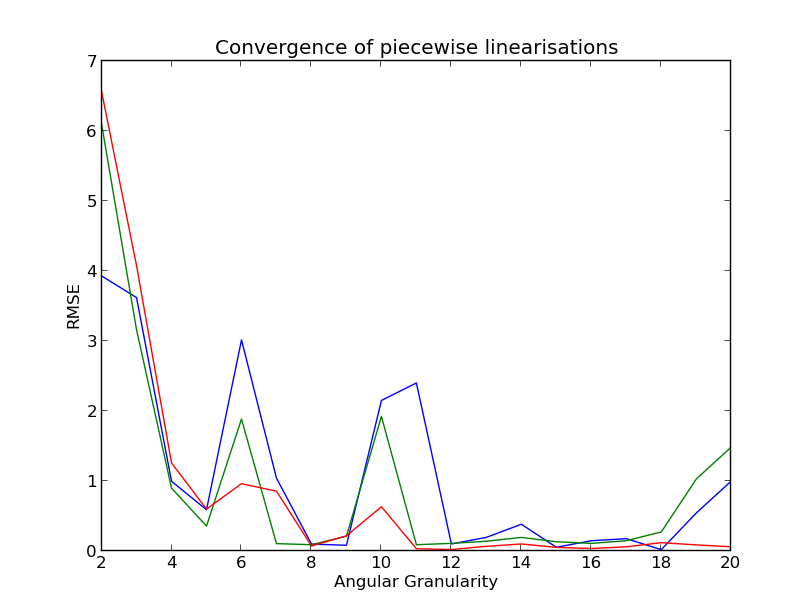}
\caption{The RMSE of voltages on Garver6y for piece-wise linearisations as a function of 
the number of segments used in the piece-wise linearisation of the voltage
angle per each 90 degrees.}
\label{fig:garver6instanceconvergence}
\end{figure}

Next, we illustrate the error in the piece-wise linearisations 
on Configuration 101 of Garver6y in Figure~\ref{fig:garver6instanceconvergence}.
On the horizontal axis, we plot the number of segments
used for the piece-wise linearisation of the voltage
angle per each 90 degrees. 
On the vertical axis, we plot the
RMSE of the solution of the corresponding piece-wise linearised instance in terms of 
voltages. 
For the remaining three dimensions, i.e. voltage
magnitude, active-, and reactive power injections,
we use piece-wise linearisations with 
uniformly 1, 2, or 3 pieces,
and obtain the blue, green, and red curves in the plot, respectively.
The evolution of RMSE over the number of segments seems disappointing.
(It should not be expected to decrease monotonically, though: 
for the example of a feasible set comprising a disk in 2D,
a particular rotation of a square yields 0 error for any objective 
function parallel to an axis, while no rotation of a pentagon yields the same.)
The evolution of run-time over the number of segments is more disappointing, still:
While the piece-wise linearisation with three segments across all dimensions takes 2.20 seconds to solve using CPLEX 12.5 with default parameters,
the piece-wise linearisations with 4, 8, 12, and 16 segments across voltage angle and three segments elsewhere take 
 3.82, 10.41, 23.84, and 96.45 seconds to solve.
Notice that this is a single configuration of a 6-bus instance,
 rather than the investment problem propers.

Finally, the performance of the lift-and-branch-and-bound procedure is somewhat promising.
Even a simple, preliminary implementation traverses a tree of 15 subproblems in 59873 seconds using SeDuMi, the SDP solver.
We envision this could be sped up much further.

\section{Conclusions}

Whereas in the operations of power systems, piece-wise linearisations may soon be replaced by 
convex relaxations
\cite{lavaei2012zero,molzahn2011,Ghaddar15}, the outlook remains less clear within investment planning. 
Although polynomial optimisation allows for global optimisation in power systems with accurate models for the physics, it remains a challenge to develop solvers that would scale to realistic instances, especially considering multiple scenarios.
Considering also that the Lavaei-Low relaxation seems difficult to extend to the investment decisions, one can hardly enumerate all the possible configurations and test them with an interior-point method, and the scalability of the piece-wise linearisations is also limited, it seems worth studying the polynomial optimisation approach in more detail.

The first results obtained with the lift-and-branch-and-bound method give some indication on how to make approaches based on polynomial optimisation applicable to investment planning in power systems, which involves both, discrete decisions and accurate models of the non-convex power flow in the constraints. If this approach proves to be scalable for larger instances, it may potentially be applied to investment planning  problems beyond power systems, where there is a combined challenge of discrete investment decisions, continuous operational decisions and non-convex system dynamics, such as in gas and water network optimisation or in traffic management.
We conjecture that the lift-and-branch-and-bound method has finite convergence for a large class of instances, although we do not prove so in this paper, and envision much further research focussed on it.

\bibliographystyle{ieeetr}
\bibliography{references-jonas}

\end{document}